# Efficient augmented block designs for unreplicated test treatments

Rahul Mukerjee

Indian Institute of Management Calcutta

Joka, Diamond Harbour Road, Kolkata 700 040, India

E-mail: rmuk0902@gmail.com

*Abstract*: Augmented block designs for unreplicated test treatments are investigated under the *A*- and *MV*-criteria with respect to control versus control, test versus test and control versus test comparisons. We derive design-independent lower bounds on these criteria over a wide class of competing designs. These bounds are useful benchmarks and the resulting expressions for efficiencies enable objective assessment of any given design under the *A*- and *MV*-criteria. It is seen that the use of BIB or PBIB designs or duals thereof often leads to very high efficiencies, which compare extremely well with algorithmic constructions. Our findings also cover the case of partial replication of some test treatments. Illustrative examples, including large-scale ones, are presented.

## 1. Introduction

Augmented block designs are widely used in agricultural experiments when the primary objective is to compare a large number of early generation varieties which represent the test treatments. Available resources for the test treatments, such as the quantity of seed, are typically quite limited, permitting the use of only one plot for each of them. Thus, the test treatments are unreplicated. In addition, the experiment includes replicated controls which enable statistically valid comparisons eliminating the block effects. Thus, one starts with a block design for the controls and then augments each block there by adding a number of test treatments, with each test treatment appearing in exactly one block. Assignment of controls and test treatments to the plots is duly randomized. As indicated above, contrasts, i.e., comparisons, between the test treatments are of primary interest, while contrasts between control and test treatments are of secondary interest and those between the controls are relatively less important. Interestingly, the starting design for the controls, which we call the primal, alone determines the precision of all three types of comparisons. Thus, the primal plays a critical role in designing the experiment and our attention will be directed to its choice.

Augmented block designs have a long history. Early references include Federer (1956, 1961) and Steel (1958). Variance formulae for estimators of the three types of contrasts stated above, when the primal is a randomized complete block design or a balanced incomplete block (BIB) design or the dual thereof, were obtained by Federer (1961), Searle (1965) and Federer and Raghavarao (1975). The authors last mentioned also studied augmented row-column designs. Herzberg and Jarrett (2007) explored, under a fixed effects model, *A*-optimal designs for estimating pairwise contrasts between the test treatments with minimum average variance. This was, however, done in the restricted class of primals that are equireplicate and have at most two distinct block intersection numbers. Indeed, their findings critically depend on a relationship that holds only for equireplicate primals.

The present paper is inspired by a recent work of Haines (2021) who made an in-depth and unified study of augmented block designs under a fixed effects model. Explicit analytical formulae for the $A$- and $MV$-criteria were given there, considering all three types of contrasts separately and allowing the primal to be any block design. Note that the $MV$-criterion captures the maximum possible variance, just as the $A$-criterion concerns the average variance. Haines (2021) also proposed algorithmic constructions of primals, with a view to making the values of these design criteria as small as possible. As noted there, other design criteria, like the $D$- or $E$-criteria, are not really appealing in the present context because of the focus on pairwise contrasts.

We work in the same framework as Haines (2021). The main contributions of the present paper, in the light of the existing literature, are as follows.

(a) We obtain lower bounds on the $A$- and $MV$-criteria over a general class of primals which need not be equireplicate or binary, or satisfy any restriction, such as the one in Herzberg and Jarrett (2008), on block intersections. These bounds are design-independent and often quite sharp. Hence they serve as useful benchmarks and tell us if seemingly attractive $A$- or $MV$-criteria values for a given design are really so in absolute terms.

(b) On the basis of the bounds, we give expressions for design efficiencies in the general class. These are more informative than the values of design criteria themselves, and enable objective assessment of any given primal, found either combinatorially or algorithmically, under these criteria.

(c) We study the efficiencies for a wide range of choices of the primal, such a BIB design, a PBIB design or its dual. In the process, we completely enumerate the 781 PBIB designs in the Clatworthy (1973) tables, as well as their duals. A vast majority of these are found to have impressive $A$-efficiencies for all three types of contrasts. Examples show that these can also have high $MV$-efficiencies. This broad collection of designs and further designs obtained therefrom add significant flexibility to the choice of the primal, enabling us to entertain almost any pre-specification for its parameters.

In order to give an idea of how (a)-(c) work, suppose there are $v = 5$ controls and the primal has to be laid out in $b = 8$ blocks, each of size $k = 3$. Since $bk/v$ is not an integer, no equireplicate primal exists. However, as detailed later in Example 1, our flexible approach yields a non-equireplicate primal which consists of the 8 blocks

$$\{1, 2, 4\}, \{1, 2, 5\}, \{1, 3, 4\}, \{1, 3, 5\}, \{1, 4, 5\}, \{2, 3, 4\}, \{2, 3, 5\}, \{2, 4, 5\},$$

and has $A$-efficiencies as high as 0.986, 0.997 and 0.994 for the control versus control, test versus test and control versus test contrasts, respectively, irrespective of the number of test treatments allocated to each block. We re-emphasize that, as noted in (a) and (b) above, these efficiencies refer to the entire class of primals with $(b, v, k) = (8, 5, 3)$, carry more information than the A-criteria values themselves and objectively evaluate of the primal that we obtain.

For further illustration, consider the large-scale example in Haines (2021) motivated by the wheat selection trials considered in Clarke and Stefanova (2011). Here 12 controls are laid out in 12 blocks, each of size 6, and 228 test treatments are equally distributed among the blocks. In this setup, following (c) above, we find a PBIB design which compares extremely well with the algorithmically constructed ones in Haines (2021) under the A-criteria and has an edge over the latter under the MV-criteria. More details appear in Section 5. This PBIB design is obtained effortlessly from Clatworthy (1973) and does not require any computation intensive procedure for its generation.

Here is an outline of the paper. Section 2 describes our setup and presents the preliminaries. In Section 3, we study the A-criteria and derive bounds thereon, which lead to the expressions for the corresponding efficiencies. These are applied in Section 4 for evaluation of BIB and PBIB designs and their duals. The situation where the specified design parameters rule out the existence of any equireplicate primal is also taken up here. Section 5 is devoted to the MV-criterion. Up to this section, the unreplicated test treatments are supposed to be equally distributed among the blocks. Section 6 extends our results to the situation where this is not the case. Partial replication, where the available resources allow two replications for some test treatments, is also discussed in this section. The concluding remarks in Section 7 indicate several directions for future work. Finally, the appendix presents some proofs as well as lists of PBIB designs and their duals which have high efficiencies.

We mention that, like Haines (2021), we have not gone beyond block designs for unreplicated test treatments. This is because of the ease in implementation and consequent wide use of such designs, a point emphasized in Haines (2021), with a reference to Piepho and Williams (2016). However, in the concluding remarks, we observe that the idea of error functions, as used later in this paper, may simplify theoretical development even in more complex settings such as row-column designs.

## 2. Preliminaries

Consider a block design $d_0$ involving $v$ controls, $bs$ test treatments and $b$ blocks, each consisting of $k + s$ plots or experimental units, such that

(a) in each block, $k$ plots are assigned to the controls, and

(b) each of the $bs$ test treatments is unreplicated, i.e., appears exactly once in the design.

The primal, called $d$, is the subdesign formed by the plots assigned to the controls. It involves $b$ blocks, each of size $k$. The design $d_0$ is an augmentation of $d$, where each of the $b$ blocks of $d$ is augmented by $s$ plots to accommodate $s$ of the $bs$ test treatments, with no test treatment appearing in more than one block.

An ordered pair representation for the test treatments, which are of primary interest in an augmented block design, will help. For $j = 1,\ldots, b$, denote the $s$ test treatments in the $j$th block by ordered pairs $jw$, $w = 1,\ldots, s$. With $Z_{jw}$ and $\rho_{jw}$ denoting, respectively, the yield and the fixed treatment effect for the test treatment $jw$, we then consider the fixed effects linear model

$$E(Z_{jw}) = \beta_j + \rho_{jw}, \quad j = 1,\ldots, b, \quad w = 1,\ldots, s, \quad (1)$$

where $\beta_j$ is the fixed effect of the $j$th block.

Turning next to the controls which are allocated to $k$ plots in each block, let $Y_{ju}$ be the yield from the $u$th of these $k$ plots in the $j$th block. Analogously to (1), we model these yields as

$$E(Y_{ju}) = \beta_j + \Sigma_{i=1}^{v} \delta(ju,i)\tau_i, \quad j = 1,\ldots, b, \quad u = 1,\ldots, k, \quad (2)$$

where $\tau_i$ is the fixed effect of the $i$th control, and $\delta(ju,i)$ is an indicator which equals 1 if $Y_{ju}$ arises from the $i$th control and 0 otherwise.

Since contrasts between controls and/or test treatments are our only objects of interest, it is immaterial whether (1) and (2) include a general mean or not. The yields modeled in (1) and (2) are assumed to be uncorrelated and homoscedastic, with a common variance, say $\sigma^2$. Write $Z$, $Y$ and $\rho$ for the $(bs)\times 1$, $(bk)\times 1$ and $(bs)\times 1$ vectors, with lexicographically arranged elements $Z_{jw}$, $Y_{ju}$ and $\rho_{jw}$, respectively. The following proposition will streamline the development of design criteria. It is anticipated from (1) and (2) from consideration of error functions which were implicit in some earlier work (e.g., Federer and Raghavarao, 1975) but perhaps never used explicitly in the present context. A proof of Proposition 1 is sketched in Appendix 1 for completeness. .

**Proposition 1**. (a) *Any linear function of Y is the best linear unbiased estimator (BLUE) of its expectation in the augmented design $d_0$ if and only if it is the BLUE of its expectation in the primal $d$.*
(b) *Any linear function of Z is the BLUE of its expectation in the augmented design $d_0$.*
(c) *Any linear function of Y and Z is the BLUE of its expectation in the augmented design $d_0$ if and only if its part involving Y alone is the BLUE of its expectation in the primal $d$.*

As Proposition 1 suggests and later developments confirm, the primal $d$ alone determines the quality of inference on contrasts of interest and hence its proper choice is crucial. The $v\times b$ incidence matrix of $d$ is given by $N = (n_{ij})$, where $n_{ij}$ is the number of times the $i$th control appears in the $j$th block. On the other hand, the dual design of $d$, say $\tilde{d}$, is represented by the $b\times v$ incidence matrix $N^T$, with the superscript $T$ denoting transpose. Let $r_1,\ldots,r_v$ be the replication numbers of the $v$ controls, $R$ = diag$(r_1,\ldots,r_v)$, and write $I_b$ be the identity matrix of order $b$. Then

$$C = R - (1/k)NN^T, \quad \tilde{C} = kI_b - N^T R^{-1} N \quad (3)$$

are the intrablock matrices of $d$ and $\tilde{d}$, respectively. Let $C^+$ and $\tilde{C}^+$ denote their respective Moore-Penrose inverses (Rao, 1973, p. 26).

## 3. *A*-criteria, lower bounds thereon and *A*-efficiencies

*3.1 Variances and A-criteria*

As mentioned earlier, we focus attention on contrasts, considering specifically three types of contrasts, namely,

*cc*: control versus control, i.e., $\tau_i - \tau_{i*}$, $i \neq i*$,

*tt*: test versus test, i.e., $\rho_{jw} - \rho_{j*w*}$, $jw \neq j*w*$, and

*ct*: control versus test, i.e., $\tau_i - \rho_{jw}$.

Because the test treatments are of main interest, contrasts of types *tt* and *ct* are of primary and secondary importance to us, and those of type *cc* are less important; cf. Haines (2021). We assume the primal *d* to be a connected design. This ensures estimability of all these contrasts.

Proposition 2 below shows the BLUEs of the three types of contrasts and the variances thereof. Parts (a), (b) of this proposition are known (Haines, 2021; Herzberg and Jarrett, 2007), but part (c) is not available in its present form. For instance, Herzberg and Jarrett (2007) gave it only for equireplicate primals with at most two distinct block intersection numbers while the counterpart of (c) in Haines (2021) is in terms of $C^+$ rather than $\tilde{C}^+$. The current form of (c) is more convenient in studying efficiency. Moreover, use of Proposition 1, which itself is not hard to obtain, is seen to simplify the proof of Proposition 2 considerably, making all its three parts quite transparent.

To present Proposition 2, let $\beta = (\beta_1, ..., \beta_b)^T$ be the vector of block effects, $e_1, ..., e_v$ be unit vectors of order *v*, and $\tilde{e}_1, ..., \tilde{e}_b$ be unit vectors of order *b*. Also, for any $i, j, i*(\neq i)$ and $j*(\neq j)$, let

$$\xi_{ij} = \tilde{e}_j - N^T R^{-1} e_i, \qquad (4)$$

and write $W_1(i, i*)$, $W_2(j, j*)$ and $W_3(i, j)$, respectively, for the BLUEs of $\tau_i - \tau_{i*}$, $\beta_j - \beta_{j*}$ and the block contrast $\xi_{ij}^T \beta$, on the basis of *Y*, i.e., arising from the primal *d*. In addition, let $T_{ci}$ be the total yield from the $r_i$ plots receiving the *i*th control.

**Proposition 2**. *The following hold for the augmented design $d_0$.*

(a) *For every $i \neq i*$, the BLUE of $\tau_i - \tau_{i*}$ is $W_1(i, i*)$, with variance $\sigma^2 V_{cc}(i, i*)$, where $V_{cc}(i, i*) = (e_i - e_{i*})^T C^+ (e_i - e_{i*})$.*

(b) *For every $jw \neq j*w*$, the BLUE of $\rho_{jw} - \rho_{j*w*}$ is (i) $Z_{jw} - Z_{jw*}$ if $j = j*$, $w \neq w*$, and (ii) $Z_{jw} - Z_{j*w*} - W_2(j, j*)$ if $j \neq j*$, with variance $2\sigma^2$ and $\sigma^2\{2 + V_{tt}(j, j*)\}$, under (i) and (ii), respectively, where $V_{tt}(j, j*) = (\tilde{e}_j - \tilde{e}_{j*})^T \tilde{C}^+ (\tilde{e}_j - \tilde{e}_{j*})$.*

(c) *For every $i$ and $jw$, the BLUE of $\tau_i - \rho_{jw}$ is $(T_{ci}/r_i) + W_3(i,j) - Z_{jw}$, with variance $\sigma^2 V_{ct}(i,j)$,*

*where $V_{ct}(i,j) = 1 + (1/r_i) + \xi_{ij}^T \tilde{C}^+ \xi_{ij}$.*

Parts (a) and (b)(i) of Proposition 2 are immediate from (1) and Proposition 1(a),(b). Similarly, parts (b)(ii) and (c) of Proposition 2 follow readily from (1), (4) and Proposition 1(c), because

$$\rho_{jw} - \rho_{j^*w^*} = \{(\beta_j + \rho_{jw}) - (\beta_{j^*} + \rho_{j^*w^*})\} - (\beta_j - \beta_{j^*}),$$

$$\tau_i - \rho_{jw} = (\tau_i + e_i^T R^{-1} N\beta) + \xi_{ij}^T \beta - (\beta_j + \rho_{jw}),$$

and the BLUEs of $(\tau_i + e_i^T R^{-1} N\beta)$ and $\xi_{ij}^T \beta$ in $d$, namely, $T_{ci}/r_i$ and $W_3(i,j)$, are uncorrelated, with respective variances $\sigma^2/r_i$ and $\sigma^2 \xi_{ij}^T \tilde{C}^+ \xi_{ij}$; cf. Dey (2010, pp. 9-12).

Proposition 2 leads to expressions for the *A*-criteria, based on average variance, separately for contrasts of types *cc*, *tt* and *ct*. These are given by

$$A_{cc} = 2\Sigma\Sigma_{1 \leq i < i^* \leq v} V_{cc}(i,i^*)/\{v(v-1)\} = 2\text{tr}(C^+)/(v-1), \tag{5}$$

$$A_{tt}(s) = 2 + 2s^2 \Sigma\Sigma_{1 \leq j < j^* \leq b} V_{tt}(j,j^*)/\{bs(bs-1)\} = 2[1 + \{s/(bs-1)\}\text{tr}(\tilde{C}^+)], \tag{6}$$

$$A_{ct} = \Sigma_{i=1}^v \Sigma_{j=1}^b \Sigma_{w=1}^s V_{ct}(i,j)/(vbs) = \Sigma_{i=1}^v \Sigma_{j=1}^b V_{ct}(i,j)/(vb)$$

$$= 1 + (1/v)\Sigma_{i=1}^v (1/r_i) + (1/b)\text{tr}(\tilde{C}^+) + (1/v)\text{tr}(R^{-1} N \tilde{C}^+ N^T R^{-1}). \tag{7}$$

Note that $A_{cc}$ and $A_{ct}$ do not involve on *s*, but $A_{tt}(s)$ depends on *s*. While (5) agrees with Haines (2021) and (6) is in the spirit of Herzberg and Jarrett (2007), the form of $A_{ct}$ in (7) is not available in the literature. It follows from (4) noting that

$$\Sigma_{i=1}^v \Sigma_{j=1}^b \xi_{ij}^T \tilde{C}^+ \xi_{ij} = \text{tr}\{\tilde{C}^+ (\Sigma_{i=1}^v \Sigma_{j=1}^b \xi_{ij} \xi_{ij}^T)\} = \text{tr}\{\tilde{C}^+ (vI_b + bN^T R^{-1} R^{-1} N)\},$$

because $\tilde{C}^+$ has vanishing row and column sums. We remark that, upon dualization, i.e., interchange of the roles of controls and blocks, (7) matches its counterpart in Haines (2021). However, (7) involves $\tilde{C}^+$ and not $C^+$, and hence it is found advantageous in the next subsection in setting a lower bound on $A_{ct}$, as the primal $d$ has a constant block size $k$ without necessarily being equireplicate.

*3.2 Lower bounds on A-criteria and A-efficiencies*

The *A*-criteria in (5)-(7) are dictated by the primal $d$ which determines both $C^+$ and $\tilde{C}^+$. Hence we now turn to the choice of $d$ so as to achieve optimality or high efficiency with respect to these criteria by minimizing their values or keeping their values small. From this perspective, given $b$, $v$ and $k$, let $D = D(b,v,k)$ denote the class of all possible connected primals, including the ones that are not binary or equireplicate. In view of (5) and (6), existing *A*-optimality results on ordinary block designs, like those on BIB or certain types of most balanced group divisible designs (Cheng, 1978) have their parallels for contrasts of types *cc* and *tt*.

Existing optimality results are, however, sparse and often unavailable for arbitrary *b*, *v* and *k*. Then, as in Haines (2021), one may find the primal *d* algorithmically, making the *A*-criteria values of interest as small as possible. Another promising approach that we adopt here is to choose *d* as a readily available block design, which shows some kind of balance or near balance even in the absence of any known optimality result, and examine the consequences on the *A*-criteria in (5)-(7). Good candidates are two-associate class PBIB designs and their duals. There is a vast pool of such designs that can cater to a wide variety of *b*, *v* and *k*, see Clatworthy (1973).

In either of the two approaches mentioned above, namely, algorithmic construction or examination of a readily available design, it is imperative to compare the resulting *A*-criteria values with design-independent lower bounds thereon over the class *D*. The ratio of any such bound to the corresponding criterion value measures design efficiency, and a high efficiency, say 0.95 or 0.99, is more assuring about design performance in the class *D* than the criterion value itself.

Theorem 1 below presents lower bounds on the *A*-criteria in (5)-(7). Here

$$L = (v-1)^2/\{b(k-1)\}, \quad \widetilde{L} = (b-1)^2/(bk-v), \qquad H = \{h/(f+1)\} + \{(v-h)/f\}, \qquad (8)$$

where *f* is the largest integer in *bk*/*v* and $h = bk - vf$.

**Theorem 1**. *For every primal d in D* [= *D*(*b*, *v*, *k*)], (a) $A_{cc} \geq A_{cc}^{\text{bound}}$, (b) $A_{tt}(s) \geq A_{tt}^{\text{bound}}(s)$, (c) $A_{ct} \geq A_{ct}^{\text{bound}}$, *where the design-independent lower bounds are given by*

$$A_{cc}^{\text{bound}} = 2L/(v-1), \qquad A_{tt}^{\text{bound}}(s) = 2[1 + \{s/(bs-1)\}\widetilde{L}],$$

$$A_{ct}^{\text{bound}} = 1 + \{(k+1)/(vk)\}H + (1/b)\widetilde{L} - (1/bk).$$

*Proof*. For every primal *d* in *D*, the arithmetic mean-harmonic mean inequality, coupled with the inequality $n_{ij}^2 \geq n_{ij}$, yields $\text{tr}(C^+) \geq L$ and $\text{tr}(\widetilde{C}^+) \geq \widetilde{L}$. Thus, (a) and (b) follow from (5) and (6).

Turning to the proof of (c), which requires more effort, from (3) observe that no eigenvalue of $\widetilde{C}$ can exceed *k*. Hence, a spectral decomposition of $\widetilde{C}^+$ shows that $\widetilde{C}^+ - (1/k)(I_b - b^{-1}1_b1_b^T)$ is nonnegative definite (nnd), where $1_b$ is the *b*x1 vector of ones. As a result,

$$\text{tr}(R^{-1}N\widetilde{C}^+N^TR^{-1}) \geq (1/k)\text{tr}\{R^{-1}N(I_b - b^{-1}1_b1_b^T)N^TR^{-1}\}$$

$$= (1/k)\text{tr}(R^{-1}NN^TR^{-1}) - (bk)^{-1}1_b^TN^TR^{-1}R^{-1}N1_b$$

$$= (1/k)\Sigma_{i=1}^v\Sigma_{j=1}^b(n_{ij}^2/r_i^2) - (v/bk)$$

$$\geq (1/k)\Sigma_{i=1}^v\Sigma_{j=1}^b(n_{ij}/r_i^2) - (v/bk) = (1/k)\Sigma_{i=1}^v(1/r_i) - (v/bk),$$

because $N1_b = (r_1,...,r_v)^T$, and therefore $1_b^T N^T R^{-1} R^{-1} N 1_b = v$. Part (c) now follows from (7) noting that $\text{tr}(\widetilde{C}^+) \geq \widetilde{L}$ and, by (8), $\Sigma_{i=1}^{v}(1/r_i) \geq H$, for any positive integers $r_1,...,r_v$ summing to $bk$. □

By (5)-(7) and Theorem 1, the efficiencies of any primal $d$, over the entire class $D$ and with respect to the $A$-criteria for the three types of contrasts are, are measured by

$$A_{cc}\text{eff} = \frac{A_{cc}^{\text{bound}}}{A_{cc}} = \frac{L}{\text{tr}(C^+)}, \quad A_{tt}\text{eff}(s) = \frac{A_{tt}^{\text{bound}}(s)}{A_{tt}(s)} = \frac{1+\{s/(bs-1)\}\widetilde{L}}{1+\{s/bs-1)\}\text{tr}(\widetilde{C}^+)}, \quad (9)$$

$$A_{ct}\text{eff} = \frac{A_{ct}^{\text{bound}}}{A_{ct}} = \frac{1+\{(k+1)/(vk)\}H + (1/b)\widetilde{L} - 1/(bk)}{1+(1/v)\Sigma_{i=1}^{v}(1/r_i) + (1/b)\text{tr}(\widetilde{C}^+) + (1/v)\text{tr}(R^{-1}N\widetilde{C}^+ N^T R^{-1})}. \quad (10)$$

Note that $s/(bs-1)$ is maximum at $s = 1$ and hence $A_{tt}\text{eff}(s) \geq A_{tt}\text{eff}(1)$. By (8) and (9), this yields a conservative measure of the $A_{tt}$-efficiency of $d$, applicable to every $s$, as

$$A_{tt}\text{eff} = A_{tt}\text{eff}(1) = \frac{1+\{(b-1)/(bk-v)\}}{1+\{1/(b-1)\}\text{tr}(\widetilde{C}^+)}. \quad (11)$$

Indeed, not only (11) but also (9) and (10) are potentially conservative because the bounds in Theorem 1 may be unattainable in $D$. Nevertheless, as seen in the next section, even this conservative approach often allows us to identify primals with $A$-efficiencies close to 1, thus showing that the lower bounds in Theorem 1 are sharp or nearly attainable even when they are not exactly attainable. The formulae

$$C^+ = (C + v^{-1}J_{vv})^{-1} - v^{-1}J_{vv}, \quad \widetilde{C}^+ = (\widetilde{C} + b^{-1}J_{bb})^{-1} - b^{-1}J_{bb},$$

which hold for any connected primal $d$ (Dey, 2010, p. 17) facilitate the computation of these efficiencies. Here $J_{vv}$ and $J_{bb}$ are matrices of ones, of orders $v$ and $b$, respectively.

## 4. $A$-efficient designs

*4.1 Case when equireplicate primals exist*

Suppose $bk/v$ (= $r$, say) is an integer, allowing the primal $d$ to be chosen as an equireplicate design with common replication number $r$. For such an equireplicate $d$, we have the simplifying features

$$\text{tr}(\widetilde{C}^+) = (r/k)\text{tr}(C^+) + (b-v)/k, \quad \text{tr}(R^{-1}N\widetilde{C}^+ N^T R^{-1}) = (v/b)\text{tr}(\widetilde{C}^+) - (b-1)/r. \quad (12)$$

The first equation in (12) is well-known; see Roy, 1958. It enables calculation of any one of $A_{cc}\text{eff}$ and $A_{tt}\text{eff}$ in (9) and (11) from the other. On the other hand, if one recalls (3), notes that $k/r = v/b$, and considers spectral decompositions of $\widetilde{C}$ and $\widetilde{C}^+$, then the second equation in (12) is obtained as

$$\text{tr}(R^{-1}N\widetilde{C}^+ N^T R^{-1}) = (1/r)\text{tr}(\widetilde{C}^+ N^T R^{-1} N) = (1/r)\text{tr}\{\widetilde{C}^+(kI_b - \widetilde{C})\}$$

$$= (1/r)\text{tr}\{k\widetilde{C}^+ - (I_b - b^{-1}1_b 1_b^T)\} = (v/b)\text{tr}(\widetilde{C}^+) - (b-1)/r.$$

It simplifies the computation of $A_{ct}$eff in (10). We emphasize that all these efficiencies continue to be in the entire class $D = D(b, v, k)$ which includes both equireplicate and non-equireplicate primals.

Continuing with integer $bk/v$ $(= r)$, we now examine the $A$-efficiencies when the primal $d$ is a BIB or a PBIB design or its dual. As contrasts of types $tt$ and $ct$ are of primary and secondary interest, and those of type $cc$ are less important, we look for designs which should ideally have

$$A_{tt}\text{eff} \geq 0.99, \qquad A_{ct}\text{eff} \geq 0.97, \qquad A_{cc}\text{eff} \geq 0.95. \qquad (13)$$

First suppose $b$, $v$ and $k$ be such that a BIB design exists and the primal $d$ is chosen as such a design. Then $b \geq v$ by Fisher's inequality, and it is not hard to see that $\text{tr}(C^+) = L$, where $L$ is given by (8); see Dey (2010, Ch. 3). As a result, by (9), $A_{cc}\text{eff} = 1$, i.e., $A_{cc}$-optimality holds. Moreover, the minimum possible values of $A_{tt}$eff and $A_{ct}$eff, obtained over the range $2 \leq r, k \leq 20$ via (9)-(12), are seen to be 0.978 and 0.959, respectively, which are quite satisfactory. Also, over $2 \leq r, k \leq 20$, in keeping with (13), we find that $A_{tt}\text{eff} \geq 0.99$ and $A_{ct}\text{eff} \geq 0.97$, whenever $k \geq 3$. Note that this range of $r$ and $k$ is sufficiently wide for most practical purposes.

Next, consider the situation where $b$, $v$ and $k$ allow choosing the primal $d$ as the dual of a BIB design. Then $v \geq b$, the primal is a linked block design and $\text{tr}(\widetilde{C}^+) = \widetilde{L}$, where $\widetilde{L}$ is given by (8). Therefore, by (11), $A_{tt}\text{eff} = 1$, i.e., $A_{tt}$-optimality holds. Furthermore, the minimum possible values of $A_{ct}$eff and $A_{cc}$eff, over $2 \leq r, k \leq 20$, turn out to be 0.976 and 0.941, respectively, so that in keeping with (13), we always get $A_{ct}\text{eff} \geq 0.97$. Also, over this range of $r$ and $k$, we find that $A_{cc}\text{eff} \geq 0.95$, except for only eight combinations of $(b, v, k)$, as given by $v = \binom{b}{2}$, $k = b - 1$, $5 \leq b \leq 12$. In all these eight cases, $d$ is the dual of a BIB design with block size only two.

As seen above, choosing the primal as a BIB design its dual leads to very impressive $A$-efficiencies. However, given $b$, $v$ and $k$, often such a choice is not possible due to the stringent conditions on existence of a BIB design. This prompts us to consider PBIB designs which are more abundant, a convenient repository being Clatworthy (1973) which exhibits 781 PBIB designs with two associate classes. A complete enumeration of these designs shows that if $b$, $v$ and $k$ allow their use as the primal $d$, then as many as 505 (i.e., 65%) of them are highly $A$-efficient in the sense of satisfying (13). Moreover, another 140 (i.e., 18%) of them do not meet (13) but satisfy the slightly less demanding condition

$$A_{tt}\text{eff} \geq 0.97, \qquad A_{ct}\text{eff} \geq 0.95, \qquad A_{cc}\text{eff} \geq 0.93. \qquad (14)$$

The rest do not meet (14) and, as with BIB designs, many of these have $k = 2$.

The picture is even better when $b$, $v$ and $k$ allow choosing duals of the PBIB designs in Clatworthy (1973) as the primal $d$. In this case, 542 (i.e., 69%) of them meet (13) and are highly $A$-efficient, while another 123 (i.e., 16%) do not satisfy (13) but meet (14) and hence have attractive $A$-efficiencies. The rest do not satisfy (14) and again many of these have $k = 2$.

In Appendix 2, we list Clatworthy (1973) designs which satisfy (13) or (14), separately for the primal $d$ chosen as (I) these designs, or (II) their duals. Table 1 illustrates the $A$-efficiencies arising from some selected BIB and PBIB designs. The designs SR75, SR91, R106 and T33 in this table are available in Clatworthy (1973), while the BIB design there is the dual of the Clatworthy design SR75. All $A$-efficiencies in Table 1 are close to 1. Due to the conservative nature of our $A$-efficiency measures as mentioned in the last section, these figures can even be indicative of actual optimality, e.g., for $(b, v, k) = (12, 16, 8)$, the primal SR91, with $A_{cc}$eff $= 0.996$, $A_{tt}$eff $= 0.999$ and $A_{ct}$eff $= 0.998$, can well be $A_{cc}$-, $A_{tt}$- and $A_{ct}$-optimal in $D(12,16, 8)$. As a passing remark, this primal performs so well despite having as many as four distinct block intersection numbers. Thus, in contrast to Herzberg and Jarrett (2007), having primals with at most two distinct block intersections is not crucial to us.

Table 1. *A-efficiencies arising from some selected BIB and PBIB designs*

| Serial no. | $(b, v, k)$ | Primal $d$ | $A_{cc}$eff | $A_{tt}$eff | $A_{ct}$eff |
|---|---|---|---|---|---|
| 1 | (30, 25, 5) | BIB design | 1.000 | 0.999 | 0.996 |
| 2 | (25, 30, 6) | SR 75 | 0.995 | 1.000 | 0.996 |
| 3 | (12, 16, 8) | SR91 | 0.996 | 0.999 | 0.998 |
| 4 | (16, 12, 6) | Dual of SR91 | 0.996 | 0.999 | 0.998 |
| 5 | (20, 10, 4) | R106 | 0.988 | 0.997 | 0.994 |
| 6 | (10, 20, 8) | Dual of R106 | 0.985 | 0.998 | 0.995 |
| 7 | (10, 10, 4) | T33 | 0.988 | 0.997 | 0.991 |

*4.2 Case when equireplicate primals do not exist*

Now suppose $bk/v$ is not an integer, so that the primal $d$ cannot be an equireplicate design. In this case, we suggest starting with a $b^*$, not far from the desired $b$, such that $b^*k/v$ is an integer, and a primal $d^*$ in $D(b^*, k, v)$ having high $A$-efficiencies. This $d^*$ can be a BIB design, a PBIB design or dual thereof, as discussed in the last subsection. Then delete $b^* - b$ blocks from $d^*$ or repeat $b - b^*$ blocks of $d^*$, according as whether $b$ is less or greater than $b^*$, respectively, to get a primal $d$ in $D(b, k, v)$. Typically, the $d$ so obtained has high $A$-efficiencies, provided the blocks deleted or repeated are chosen judiciously. For example, if $|b^* - b| \geq 2$, then these deleted or repeated blocks should not have much overlap among themselves, so as to make $d$ nearly equireplicate.

Examples 1-3 and many others not reported here illustrate that the above approach works quite well. In all these examples, the $A$-efficiencies of $d^*$ refer to the class $D(b^*, k, v)$, while those of any primal obtained from $d^*$ via deletion or repetition of blocks refers to $D(b, k, v)$ for the corresponding $b$. Incidentally, one may wonder if a complete enumeration of all possibilities for deletion or repetition can yield even better results. Our computational experience suggests that the resulting gains, if any, over the idea of keeping the overlap small is nominal.

**Example 1.** Let $v = 5$, $k = 3$, and consider five possibilities for $b$, namely, $b = 8, 9, 10, 11, 12$. Then $bk/v$ is an integer only for $b = 10$. Take $b^* = 10$, and start with the primal $d^*$ as given by the BIB design with ten blocks

{1, 2, 3}, {1, 2, 4}, {1, 2, 5}, {1, 3, 4}, {1, 3, 5}, {1, 4, 5}, {2, 3, 4}, {2, 3, 5}, {2, 4, 5}, {3, 4, 5}.

This $d^*$ has $A_{cc}\text{eff} = 1$, $A_{tt}\text{eff} = 0.998$ and $A_{ct}\text{eff} = 0.994$. Its first and last blocks overlap in only one treatment which is minimum possible. Accordingly, from $d^*$, obtain primals for $b = 8, 9, 11$ and $12$ by (i) deleting the first and last blocks, (ii) deleting only the first block, (iii) repeating only the first block, and (iv) repeating the first and last blocks. For the primals so obtained, the triplet ($A_{cc}\text{eff}$, $A_{tt}\text{eff}$, $A_{ct}\text{eff}$) equals (i) (0.986, 0.997, 0.994), (ii) (0.988, 0.997, 0.994), (iii) (0.992, 0.998, 0.995) and (iv) (0.994, 0.998, 0.995). Thus, these non-equireplicate primals well inherit the high $A$-efficiencies of $d^*$ and, in fact, the $A_{ct}$-efficiencies in (iii) and (iv) slightly surpass that of $d^*$ itself. □

**Example 2.** Let $v = 30$, $k = 6$, and consider $b = 23, 24, 25, 26, 27$. Then $bk/v$ is an integer only for $b = 25$. Take $b^* = 25$ and $d^*$ as the Clatworthy (1973) design SR75. As seen in Table 1, it entails $A_{cc}\text{eff} = 0.995$, $A_{tt}\text{eff} = 1$, $A_{ct}\text{eff} = 0.996$. This $d^*$ is the dual of the BIB design in Table 1 and any two blocks of $d^*$ have one common treatment. In order to obtain primals for $b = 23, 24, 26$ and $27$ from $d^*$, (i) delete the first and second blocks, (ii) delete only the first block, (iii) repeat only the first block, and (iv) repeat the first and second blocks. As in the last example, the resulting primals have high $A$-efficiencies, with ($A_{cc}\text{eff}$, $A_{tt}\text{eff}$, $A_{ct}\text{eff}$) given by (i) (0.974, 0.999, 0.994), (ii) (0.985, 1, 0.996), (iii) (0.989, 0.999, 0.996) and (iv) (0.984, 0.999, 0.995). □

**Example 3.** Let $v = 10$, $k = 4$, and consider $b = 8, 9, 10, 11, 12$. Then $bk/v$ is an integer only for $b = 10$. Take $b^* = 10$ and $d^*$ as the Clatworthy (1973) design T33. As seen in Table 1, it leads to $A_{cc}\text{eff} = 0.988$, $A_{tt}\text{eff} = 0.997$ and $A_{ct}\text{eff} = 0.991$. Unlike in the last two examples, $d^*$ is neither a BIB design nor the dual of a BIB design. Because the first and last blocks of $d^*$ have one treatment in common which is minimum possible, primals for $b = 8, 9, 11, 12$ can be obtained from $d^*$ as in Example 1. For these primals, the triplet ($A_{cc}\text{eff}$, $A_{tt}\text{eff}$, $A_{ct}\text{eff}$) equals (i) (0.930, 0.995, 0.978), (ii) (0.960, 0.996, 0.988), (iii) (0.973, 0.996, 0.990) and (iv) (0.966, 0.996, 0.998). All these, except the $A_{cc}\text{eff}$ in (i), are impressive. Even the $A_{cc}\text{eff}$ in (i) meets (14). □

## 5. MV-efficiency

We now turn to the $MV$-criteria, based on maximum possible variance, separately for contrasts of types $cc$, $tt$ and $ct$. By Proposition 2, these are given by

$$MV_{cc} = \max\{V_{cc}(i, i^*) : 1 \leq i < i^* \leq v\}, \qquad MV_{tt} = 2 + \max\{V_{tt}(j, j^*) : 1 \leq j < j^* \leq b\},$$

$$MV_{ct} = \max\{V_{ct}(i, j) : 1 \leq i \leq v, 1 \leq j \leq b\}.$$

While existing *MV*-optimality results on ordinary block designs, such as BIB designs or certain types of group divisible designs (Takeuchi, 1961) have their parallels for $MV_{cc}$ and $MV_{tt}$, these results are sparse. So, for general $b$, $v$ and $k$, as in the last two sections, we consider design-independent lower bounds and resulting *MV*-efficiencies which help in assessing the performance of any primal and hence that of the associated augmented block design. These *MV*-efficiencies refer to the entire class $D = D(b, v, k)$ which includes equireplicate as well as non-equireplicate primals.

A hurdle now is that even for ordinary block designs, lower bounds specifically for the *MV*-criterion are very rare. Nonetheless, because maximum can never be less than average, it follows that $MV_{cc}$, $MV_{tt}$ and $MV_{ct}$ are bounded below, respectively, by $A_{cc}^{bound}$, $A_{tt}^{bound}(1)$ and $A_{ct}^{bound}$, as defined in Theorem 1. For instance, recalling that $\text{tr}(\tilde{C}^+) \geq \tilde{L}$, we have analogously to (5),

$$MV_{tt} \geq 2 + 2\Sigma\Sigma_{1 \leq j < j^* \leq b} V_{tt}(j, j^*)/\{b(b-1)\}$$

$$= 2[1 + \{1/(b-1)\}\text{tr}(\tilde{C}^+)] \geq 2[1 + \{1/(b-1)\}\tilde{L}],$$

where the term on the extreme right equals $A_{tt}^{bound}(1)$. In view of the above, the efficiencies of any primal $d$, over the entire class $D$ and with respect to the *MV*-criteria for the three types of contrasts, can be measured by

$$MV_{cc}\text{eff} = \frac{A_{cc}^{bound}}{MV_{cc}}, \qquad MV_{tt}\text{eff} = \frac{A_{tt}^{bound}(1)}{MV_{tt}}, \qquad MV_{ct}\text{eff} = \frac{A_{ct}^{bound}}{MV_{ct}}. \qquad (15)$$

Although (15) appears to be rather conservative, the next two examples show that it can allow us to identify designs with high *MV*-efficiencies.

**Example 4**. For the designs serially numbered 1,..., 7 in Table 1, the triplet ($MV_{cc}$eff, $MV_{tt}$eff, $MV_{ct}$eff) equals

(1.000, 0.994, 0.984),  (0.967, 1.000, 0.985),  (0.938, 0.995, 0.981),  (0.962, 0.990, 0.980),
(0.900, 0.984, 0.969),  (0.930, 0.986, 0.971),  (0.964, 0.992, 0.952),

respectively. Thus, all of them have $MV_{tt}$eff > 0.98 and $MV_{ct}$eff > 0.95. Moreover, $MV_{cc}$eff, which is of lesser interest, exceeds 0.95 for four of them and is not less than 0.90 for any of them. □

**Example 5**. We consider the large-scale example in Haines (2021, Subsection 4.2.1), where $b = v = 12$, $k = 6$, $s = 19$, i.e., there are $bs = 228$ test treatments. By algorithmic construction, with a view to keeping the $A_{tt}$- and $A_{ct}$-criteria values as small as possible in this situation, Haines (2021) obtained two primals $d_1$ and $d_2$, of which the first one is non-equireplicate and the second one is equireplicate. On the other hand, from Clatworthy (1973) tables, we find four PBIB designs, namely, S29, SR68 and two solutions for SR67, which may be considered as primals for $(b, v, k) = (12, 12, 6)$. Of these, SR68 performs better than the others under the *A*- and *MV*-criteria. Table 2 compares SR68 with the algorithmically obtained primals $d_1$ and $d_2$. In this table, we report $A_{tt}(s)$ and $A_{tt}\text{eff}(s)$ for the specific $s$ (=19) considered here. Each efficiency is shown in parentheses below the corresponding criterion value and two typos in Haines (2021), in $A_{ct}$ for $d_1$ and $MV_{tt}$ for $d_2$, are corrected. Table 2 shows that SR68 has the same $A_{tt}(s)$ as $d_1$ and the same $A_{ct}$ as $d_2$, thus equaling these two designs under the respective criteria that their algorithmic construction was directed to. In general, the primal SR68, obtained effortlessly from the Clatworthy (1973) tables, always compares very well with the computation intensively generated $d_1$ and $d_2$, and even has an edge over the latter under all the $MV$-criteria. Moreover, our design-independent lower bounds and resulting $A$- and $MV$-efficiencies objectively assure us about the performance of the primals, a feature is not shared by empirical stopping rules based on these criteria values alone in algorithmic construction. □

Table 2. *A- and MV-criteria values and efficiencies for* $(b, v, k, s) = (12, 12, 6, 19)$

| Primal | $A_{cc}$ | $A_{tt}(s)$ | $A_{ct}$ | $MV_{cc}$ | $MV_{tt}$ | $MV_{ct}$ |
|---|---|---|---|---|---|---|
| $d_1$ | 0.388 | 2.338 | 1.361 | 0.571 | 2.378 | 1.482 |
|  | (0.944) | (0.999) | (0.991) | (0.642) | (0.995) | (0.910) |
| $d_2$ | 0.368 | 2.339 | 1.351 | 0.377 | 2.388 | 1.387 |
|  | (0.996) | (0.999) | (0.998) | (0.972) | (0.991) | (0.972) |
| SR68 | 0.367 | 2.338 | 1.351 | 0.375 | 2.375 | 1.382 |
|  | (0.998) | (0.999) | (0.998) | (0.978) | (0.996) | (0.979) |

The last two examples suggest that the $MV$-efficiency measures in (15) are useful benchmarks when $bk/v$ is an integer, which is the case in these examples. If $bk/v$ is not an integer, then the outcome is of mixed nature. As an illustration, if $v = 5$ and $k = 3$, then for $b = 8, 9, 11$ and 12, the primals obtained in Example 1 have ($MV_{cc}\text{eff}, MV_{tt}\text{eff}, MV_{ct}\text{eff}$) equal to (0.903, 0.983, 0.903), (0.889, 0.984, 0.925), (0.909, 0.986, 0.936) and (0.941, 0.985, 0.946), respectively. While $MV_{tt}\text{eff} > 0.98$ in all these, the $MV_{cc}\text{eff}$ and $MV_{ct}\text{eff}$ values are all less than 0.95 and some of them are even close to 0.90. Thus, either (i) these designs do not really perform so well under the $MV_{cc}$- and $MV_{ct}$-criteria, or (ii) the $MV_{cc}$- and $MV_{ct}$-efficiency measures in (15) are themselves too conservative when $bk/v$ is not an integer. Given the high $A_{cc}$- and $A_{ct}$-efficiencies that these designs were seen to enjoy in Example 1, possibility (i) seems unlikely. This underscores the need for developing sharper bounds on the $MV$-criteria and hence less conservative $MV$-efficiency measures, when $bk/v$ is not an integer.

**6. Two extensions**

*6.1 Unequally distributed test treatments*

So far, we considered the situation where the test treatments are equally distributed among the blocks. This may not always be possible, e.g., the number of test treatments may not be an integral multiple of the number, $b$, of blocks. More generally, now suppose there are $S = s_1 + ... + s_b$ test treatments of which $s_j$ appear in the $j$th block ($j = 1,\ldots, b$), all test treatments being unreplicated. In this case, it is evident from Proposition 2 that the $MV$-criteria remain the same as in the last section and hence these is no need to study these afresh. The same holds for the $A_{cc}$-criterion as well.

However, the $A_{tt}$- and $A_{ct}$-criteria get affected. Recalling Proposition 2, the expressions for $A_{tt}(s)$ and $A_{ct}$ in (6) and (7) now change to

$$A_{tt}(s^\#) = 2 + 2\Sigma\Sigma_{1\le j< j^*\le b}\, s_j s_{j^*} V_{tt}(j, j^*)/\{S(S-1)\}, \qquad (16)$$

$$A_{ct}(s^\#) = \Sigma_{i=1}^{v}\Sigma_{j=1}^{b}\Sigma_{w=1}^{s_j} V_{ct}(i,j)/(vS) = \Sigma_{i=1}^{v}\Sigma_{j=1}^{b} s_j V_{ct}(i,j)/(vS), \qquad (17)$$

where $s^\#$ stands for $s_1,\ldots,s_b$ collectively. As shown in Appendix 3, these have design-independent lower bounds as given, respectively, by

$$A_{tt}^{\text{bound}}(s^\#) = 2 + \frac{(4/k)\Sigma\Sigma_{1\le j<j^*\le b}\phi_{jj^*} + 2s_0^2 b\tilde{L}}{b\bar{s}(b\bar{s}-1)}, \qquad (18)$$

$$A_{ct}^{\text{bound}}(s^\#) = 1 + \frac{(k\bar{s}+s_0)H}{vk\bar{s}} + \frac{s_0\tilde{L}}{b\bar{s}} - \frac{s_0}{bk\bar{s}}, \qquad (19)$$

where $s_0 = \min(s_1,\ldots,s_b)$, $\bar{s} = (s_1 + \ldots + s_b)/b$ and $\phi_{jj^*} = s_j s_{j^*} - s_0^2$ ($j, j^* = 1,\ldots,b$). These bounds reduce to their counterparts in Theorem 1 when $s_1 = \ldots = s_b = s$.

For any given $s_1,\ldots,s_b$, in view of the above, the $A_{tt}$- and $A_{ct}$-efficiencies of any primal $d$, over the class $D(b, v, k)$ can now be measured by

$$A_{tt}\text{eff}(s^\#) = \frac{A_{tt}^{\text{bound}}(s^\#)}{A_{tt}(s^\#)}, \qquad A_{ct}\text{eff}(s^\#) = \frac{A_{ct}^{\text{bound}}(s^\#)}{A_{ct}(s^\#)}.$$

These efficiency measures are again conservative. However, as seen below, they can be close to 1 and thus facilitate identification of highly efficient designs.

We now revisit Example 5 and examine how the primals $d_1$, $d_2$ and SR68 there perform when the test treatments are unequally distributed among the blocks. This is done for $S = 100$, 150 and 200, none of which is an integral multiple of $b = 12$. In each case, we take $s_1,\ldots,s_b$ as nearly equal as possible, e.g., for $S = 100$, we consider $s_1 = \ldots = s_8 = 8$, $s_9 = \ldots = s_{12} = 9$. Interestingly, it turns out that for all three choices of $S$, the pair ($A_{tt}\text{eff}(s^\#)$, $A_{ct}\text{eff}(s^\#)$) equals (0.999, 0.986) for $d_1$, (0.998, 0.993) for $d_2$, and (0.999, 0.993) for SR68. Thus, the primal SR68 from Clatworthy (1973) again compares very well with the algorithmically obtained $d_1$ and $d_2$. Moreover, as seen in Table 2, it continues to have an edge over the other two under the $A_{cc}$- and $MV$-criteria, which do not depend on whether $s_1,...,s_b$ are equal or not. Further computations, not reported here, indicate that the same pattern persists even when $s_1,...,s_b$ vary among themselves more substantially.

*6.2 Partially replicated test treatments*

Hitherto in this paper, all test treatments were unreplicated. Now, suppose available resources allow two replications for some of them, and controls are not used. Then there are three types of treatment contrasts, *rr*, *rt* and *tt*, all of primary interest, where *r* and *t* refer, respectively, to the twice-replicated and unreplicated test treatments. Expressions for the corresponding *A*- and *MV*-criteria and efficiencies can be obtained as in the preceding sections, with controls formally replaced by the twice-replicated test treatments. We skip the details to avoid repetition but visit Haines (2021, Subsection 4.2.2) to see how use of known block deigns can again compete very well with algorithmic constructions.

Following Haines (2021), in the setup of Example 5, suppose there are no controls, but consider 36 twice-replicated test treatments, in addition to the 228 unreplicated ones. Each of the $b=12$ blocks now includes $s = 19$ of the unreplicated test treatments and has 6 additional plots for the twice-replicated test treatments. These additional plots form a subdesign, say $d_{\text{rep}}$, involving the 36 twice-replicated test treatments and 12 blocks, each of size 6. The subdesign $d_{\text{rep}}$, analogous to the primal *d* for controls in the previous sections, dictates the *A*- and *MV*-criteria. Figure 2(b) in Haines (2021) shows the algorithmically obtained best $d_{\text{rep}}$ under the $A_{rr}$-criterion. We find that this $d_{\text{rep}}$ is isomorphic to the design LS74 in Clatworthy (1973), and hence the latter is as good as the former under all the design criteria. Thus, even with twice-replicated test treatments, our approach of examining known designs can help in averting computational construction.

**7. Concluding remarks**

Augmented block designs for unreplicated test treatments were investigated in this paper under the *A*- and *MV*-criteria. Design-independent lower bounds on these criteria and resulting expressions for efficiencies, as derived here, enabled us to assess objectively the performance of these designs. It was also seen that choosing the primal as a BIB or a PBIB design or the dual thereof often leads to very high efficiencies, and competes very well with algorithmic constructions. Several directions for future research emerge from the present work.

The case when $bk/v$ is not an integer deserves further attention. Then no equireplicate primal exists and it will be worthwhile to study if any direct construction performs even better than addition or deletion of blocks as done here. Moreover, as noted in Section 5, there is a need to obtain sharper lower bounds on the *MV*-criteria when $bk/v$ is not an integer. This is challenging even for ordinary block designs but one may examine if, for example, the approaches in Jacroux (1983) or Roy and Shah (1984) can be extended to obtain sufficiently general results.

The idea of error functions, which is inherent in Proposition 1 and holds for linear models in general, should also facilitate a theoretical study of augmented row-column designs for unreplicated test treatments. It will be of interest to investigate if, as in the present work, this can lead to expressions for the $A$- and $MV$-criteria for such designs and lower bounds thereon. These design criteria are again expected to depend only on the primal, i.e., the row-column subdesign for controls with cells left empty for the test treatments. The associated bounds can serve, like here, as useful benchmarks for objective evaluation of any such primal. For example, following a proposal in Haines (2021), the performance of space-filling designs as primals may be studied. One may as well wish to explore how well the design-criteria values for algorithmically generated designs, such as those in Vo-Thanh and Piepho (2020), compare with the corresponding lower bounds, if available as above.

It should be possible to exploit error functions to obtain analogs of the present results under a fixed effects model, with an intraclass correlation structure for yields from the same block. Further work on how far this approach or its modifications apply to more challenging correlation structures (see Cullis et al., 2020, and the references there) will be welcome.

We conclude with the hope that the present work will generate interest in the above and related problems.

**Appendix 1: Proof of Proposition 1**

An error function is a linear function of $Y$ and $Z$, say, $l_1^T Y + l_2^T Z$, with expectation identically equal to zero. By (1) and (2), this expectation equals $l_2^T \rho$ plus a linear function of the block and control effects. Thus, $l_1^T Y + l_2^T Z$ is an error function if and only if (i) $l_2 = 0$ and (ii) $E(l_1^T Y)$ is identically equal to zero. So, the only error functions in the augmented design $d_0$ are the ones that involve $Y$ alone and are error functions in the primal $d$. The result now follows because any linear observational function is the BLUE its expectation if and only if it is uncorrelated with every error function (Rao, 1973, p. 308). In parts (b), (c), we also use the fact that $Y$ and $Z$ are uncorrelated. □

**Appendix 2: $A$-efficient PBIB designs**

The Clatworthy (1973) tables exhibit eight types of PBIB designs, each with two associate classes, namely, singular (S), semi-regular (SR), regular (R), triangular (T), Latin square type (LS), cyclic (C), partial geometries (PG) and miscellaneous (M). We consider choosing the primal $d$ as (I) these PBIB designs or (II) their duals, and separately for (I) and (II), list PBIB designs which lead to $A$-efficient augmented block designs as a result.

(I) *Primal d chosen as a PBIB design*

(a) The following 505 PBIB designs satisfy $A_{tt}\text{eff} \geq 0.99$, $A_{ct}\text{eff} \geq 0.97$ and $A_{cc}\text{eff} \geq 0.95$.

```
S:   1-5, 7, 8, 18-34, 36-41, 44, 46, 49-72, 74-76, 80-121, 124
SR:  18, 20-25, 27, 29, 32, 34-110
```

R: 16, 17, 42-48, 50-52, 54-59, 62, 65, 67, 68, 70-72, 74, 75, 77-79, 81-84, 87-103, 106-114, 116-120, 122, 124-142, 144, 145, 147-150, 151, 153-157, 159-167, 169-177, 179-187, 189-194, 196, 197, 199-205, 208, 209
T: 12, 13, 18, 33-37, 40, 41, 43-47, 54, 56-64, 70, 71, 77, 83-85, 91-95
LS: 12-15, 26, 27, 29, 31, 33, 36-44, 48-50, 52, 53, 55, 56, 58, 61-66, 68, 69, 71-73, 75, 77, 79, 82, 83, 85-88, 90-92, 94, 97-101, 103-107, 109-111, 113, 116-118, 120-125, 127-129, 131, 134-136, 138, 139, 141, 142, 144
C: 12-15, 19, 21-29
PG: 5-15, 6a, 13a
M: 16-18, 20-22, 25-42, 31a, 32a, 37a

(b) The following 140 PBIB designs do not meet at least one of the conditions $A_{tt}\text{eff} \geq 0.99$, $A_{ct}\text{eff} \geq 0.97$, $A_{cc}\text{eff} \geq 0.95$, but satisfy $A_{tt}\text{eff} \geq 0.97$, $A_{ct}\text{eff} \geq 0.95$ and $A_{cc}\text{eff} \geq 0.93$.

S: 6, 9-12, 35, 42, 43, 45, 48, 73, 77-79, 122, 123
SR: 8, 10, 12, 16, 17, 19, 26, 28, 30, 31, 33
R: 1, 3, 6, 7, 9, 10, 12, 13, 19-21, 23-25, 27, 30-32, 35-37, 39-41, 60, 69, 73, 76, 80, 85, 86, 104, 105, 158, 195, 206
T: 9-11, 14, 15, 20-32, 38, 39, 42, 48-53, 55, 65-69, 72-76, 78-82, 86-90, 96-100
LS: 17, 20, 22, 47, 137, 140, 143, 145, 146
C: 6, 7, 9, 20
PG: 2-4
M: 3, 13-15, 15a, 19, 23, 24

(II) *Primal d chosen as the dual of a PBIB design*

(a) The following 542 PBIB designs satisfy $A_{tt}\text{eff} \geq 0.99$, $A_{ct}\text{eff} \geq 0.97$ and $A_{cc}\text{eff} \geq 0.95$.

S: 1-5, 7, 8, 10, 12, 18-34, 36, 37, 39-41, 44, 46, 48-68, 70-72, 74, 76, 80-111, 113, 115, 118, 119, 121, 124
SR: 19-110
R: 3, 6, 7, 9, 10, 13, 16, 17, 24, 27, 42-52, 54-59, 62, 65, 67, 68, 70-79, 81-103, 105-114, 116-122, 124-142, 144-151, 153-167, 169-177, 179-187, 189-194, 196, 197, 199-205, 208, 209
T: 11-13, 18, 21, 22, 30-38, 40-47, 51-64, 69-71, 76, 77, 82-85, 90-95
LS: 12-15, 17, 20, 26, 27, 31, 33, 37-44, 48-50, 53, 55, 56, 58, 61-66, 68, 69, 71-73, 77, 79, 82, 83, 86-88, 90-92, 94, 97-101, 103-107, 109-111, 113, 116-118, 120-125, 127-129, 131, 134-136, 138, 139, 141, 142, 144
C: 9, 12-15, 19-29
PG: 5-15, 6a, 13a
M: 13-15, 15a, 17-22, 24-42, 31a, 32a, 37a

(b) The following 123 PBIB designs do not meet at least one of the conditions $A_{tt}\text{eff} \geq 0.99$, $A_{ct}\text{eff} \geq 0.97$, $A_{cc}\text{eff} \geq 0.95$, but satisfy $A_{tt}\text{eff} \geq 0.97$, $A_{ct}\text{eff} \geq 0.95$ and $A_{cc}\text{eff} \geq 0.93$.

S: 6, 9, 11, 13-17, 35, 38, 42, 43, 45, 47, 69, 73, 75, 77-79, 112, 114, 116, 117, 120, 122, 123
SR: 3-5, 8, 12, 18
R: 1, 12, 15, 19, 20, 23, 25, 30-33, 35-37, 39-41, 53, 60, 61, 69, 80, 104, 115, 123, 195, 206
T: 9, 10, 15, 19, 20, 23-29, 39, 48-50, 65, 67, 68, 72-75, 78, 80, 81, 86, 88, 89, 96, 99, 100
LS: 10, 11, 22, 23, 29, 35, 36, 46, 47, 52, 60, 67, 70, 75, 81, 85, 96, 115, 137, 146
C: 6-8, 17, 18
PG: 2-4
M: 9, 16, 23

**Appendix 3: Lower bounds on $A_{tt}$- and $A_{tc}$-criteria for unequally distributed test treatments**

(a) About $A_{tt}^{\text{bound}}(s^{\#})$: As in the proof of Theorem 1, $\text{tr}(\tilde{C}^{+}) \geq \tilde{L}$ and $\tilde{C}^{+} - (1/k)(I_b - b^{-1}1_b 1_b^T)$ is nnd, so that $V_{tt}(j, j^*) \geq 2/k$. Also, analogously to (6), $\Sigma\Sigma_{1 \leq j < j^* \leq b} V_{tt}(j, j^*) = b\text{tr}(\tilde{C}^{+})$. Hence,

$$\Sigma\Sigma_{1 \leq j < j^* \leq b} s_j s_{j^*} V_{tt}(j, j^*) = \Sigma\Sigma_{1 \leq j < j^* \leq b} \phi_{jj^*} V_{tt}(j, j^*) + s_0^2 \Sigma\Sigma_{1 \leq j < j^* \leq b} V_{tt}(j, j^*)$$

$$= \Sigma\Sigma_{1 \leq j < j^* \leq b} \phi_{jj^*} V_{tt}(j, j^*) + s_0^2 b\text{tr}(\tilde{C}) \geq (2/k)\Sigma\Sigma_{1 \leq j < j^* \leq b} \phi_{jj^*} + s_0^2 b\tilde{L}.$$

If we use the above in (16) and note that $S = b\bar{s}$, then it follows that $A_{tt}(s^{\#})$ is bounded below by $A_{tt}^{\text{bound}}(s^{\#})$ in (18).

(b) About $A_{ct}^{\text{bound}}(s^{\#})$: Since $S = b\bar{s}$, if we recall the form of $V_{ct}(i, j)$ in Proposition 2, then from (17), we obtain

$$A_{ct}(s^{\#}) = \Sigma_{i=1}^{v}\Sigma_{j=1}^{b} s_j \{1 + (1/r_i) + \xi_{ij}^T \tilde{C}^{+} \xi_{ij}\}/(vb\bar{s})$$

$$= 1 + (1/v)\Sigma_{i=1}^{v}(1/r_i) + \Sigma_{i=1}^{v}\Sigma_{j=1}^{b} s_j \xi_{ij}^T \tilde{C}^{+} \xi_{ij}/(vb\bar{s})$$

$$\geq 1 + (1/v)\Sigma_{i=1}^{v}(1/r_i) + s_0 \Sigma_{i=1}^{v}\Sigma_{j=1}^{b} \xi_{ij}^T \tilde{C}^{+} \xi_{ij}/(vb\bar{s})$$

$$= 1 + (1/v)\Sigma_{i=1}^{v}(1/r_i) + \{s_0/(b\bar{s})\}\text{tr}(\tilde{C}^{+}) + \{s_0/(v\bar{s})\}\text{tr}(R^{-1}N\tilde{C}^{+}N^T R^{-1}), \quad (A.1)$$

arguing as in (7) in the last step. We again note that $\text{tr}(\tilde{C}^{+}) \geq \tilde{L}$. Furthermore, as seen while proving Theorem 1, $\Sigma_{i=1}^{v}(1/r_i) \geq H$, and

$$\text{tr}(R^{-1}N\tilde{C}^{+}N^T R^{-1}) \geq (1/k)\Sigma_{i=1}^{v}(1/r_i) - (v/bk).$$

Using the above in (A.1), it follows that $A_{ct}(s^{\#})$ is bounded below by $A_{ct}^{\text{bound}}(s^{\#})$ in (19).

**Acknowledgement**. This work was supported by a grant from the Science and Engineering Board, Government of India.